\documentclass[12pt]{amsart}                
\usepackage{amsthm}
\usepackage{amsmath}
\usepackage{amssymb}



\newtheorem*{remark}{Remark}

\newtheorem*{theorem1}{Theorem 1}
\newtheorem*{theorem2}{Theorem 2}
\newtheorem*{theorem3}{Theorem 3}
\newtheorem*{lemma1}{Lemma 1}
\newtheorem*{lemma2}{Lemma 2}
\newtheorem*{lemma3}{Lemma 3}
\newtheorem*{lemma4}{Lemma 4}

\begin{document}

\title[Conjugacy to substitution systems]{Topological conjugacy to  given constant length substitution minimal systems}

\baselineskip=17pt

\author[E. M. Coven]{Ethan M. Coven}
\author[A. Dykstra]{Andrew Dykstra}
\author[M. Keane]{Michael Keane}
\author[M. LeMasurier]{Michelle LeMasurier}

\address[Ethan M. Coven and Michael Keane]
{Department of Mathematics\\
				 Wesleyan University\\
				 Middletown CT 06459}
\email{ecoven@wesleyan.edu}
\email{mkeane@wesleyan.edu}

\address[Andrew Dykstra and Michelle LeMasurier]
				{Department of Mathematics\\
				 Hamilton College\\
				 Clinton NY 13323}
\email{adykstra@hamilton.edu}
\email{mlemasur@hamilton.edu}

\date{January 23, 2013}

\begin{abstract}
We find necessary and sufficient conditions for a symbolic dynamical system to be topologically conjugate to any given constant length substitution minimal system, thus extending the results in \cite{CKL} for the Morse and Toeplitz substitutions.
\end{abstract}

\maketitle{}

\bigskip

\section*{1. Introduction}

In \cite{CKL} three of the authors characterized 
those symbolic minimal systems that are topologically conjugate to 
the Morse Minimal System (the closure of the  shift-orbit  of the %famous 
Morse-Thue sequence) and 
those topologically conjugate to
the closely related Toeplitz Minimal System.  
The Morse result is that a symbolic minimal system
$(Y,\sigma)$ is topologically conjugate to the Morse system if and only if there exist $N \ge 1$ and $2^N$-blocks $C_0 \ne C_1$ 
such that every point in ~$Y$ can be written as a concatenation of $C_0$'s and $C_1$'s in exactly one way and such that the sequences of $C$'s in these concatenations  in some sense ``mirror''  points in the Morse system.
However the arguments in \cite{CKL} are valid only for substitutions sharing some of the properties of the Morse or Toeplitz substitutions.

Our main result is Theorem~1, which extends the Morse and Toeplitz results to all infinite substitution minimal systems generated by primitive, constant length substitutions.
In Theorem~2 we show that the ``mirroring'' property in the Morse   Dynamical Characterization Theorem in \cite{CKL} holds for the class of primitive, one-to-one substitutions $\theta$ of constant length   
at least three
such that for all $s \ne t$, $\theta(s)$ and $\theta(t)$ disagree at some place other than the first or last.
Using ideas from the proof of Theorem~1 
we show in Theorem~3 that if $\theta$ is a primitive substitution of constant length that generates an infinite system $(X_\theta,\sigma)$ and if $\zeta$ is a primitive, one-to-one  substitution of constant length that generates an infinite system $(X_\zeta,\sigma)$ topologically conjugate to $(X_\theta,\sigma)$, then the number of letters in the alphabet of ~$\zeta$ is bounded above by the number of 3-blocks that appear in
~$X_\theta$.
We give an example to show that the bound is attained for the Morse system.

\section*{2. Background}

In this paper a {\bf dynamical system} is a pair
$(X,T)$, where $X$ is a compact metric space and
$T:X \to X$ is a  
homeomorphism.
The notion of ``sameness'' for dynamical systems is
{\bf topological conjugacy}: $(X,T)$ and $(Y,S)$
are topologically conjugate if and only if
there exists a  
homeomorphism
$F: (X,T) \to (Y,S)$ such that
$F \circ T \equiv S \circ F$.
In this case $F$ is called a topological conjugacy
and $F^{-1}: (Y,S) \to (X,T)$
is also a topological conjugacy.
A {\bf topological semi-conjugacy}, also called 
a {\bf factor map}, is a continuous, onto  map
$F : (X,T) \to (Y,S)$ such that $F \circ T \equiv S \circ F$.

A dynamical system $(X,T)$ is called {\bf minimal}
if and only if it contains no proper subsystem.
Equivalently, if $X' \subseteq X$ is nonempty, closed, and $T$-invariant (i.e. $T(X') \subseteq X')$, then $X' = X$. 
Equivalently again, every orbit 
$\{T^n(x): -\infty <n < \infty\}$ is dense in $X$.

A {\bf symbolic system} is a subsystem of some
$(\prod_{-\infty}^\infty \mathcal S,\sigma)$, 
where $\mathcal S$ is a finite set of symbols, also called the 
{\bf alphabet}, and 
$\sigma :\prod_{-\infty}^\infty \mathcal S \to \prod_{-\infty}^\infty \mathcal S$
is the (left) {\bf shift}: 
$(\sigma(x))_i = x_{i+1}, -\infty < i < \infty$.
As is customary, we shall abuse notation
and write $(X,\sigma)$ instead of $(X,\sigma|_X)$.
A basic fact about symbolic systems is the
Curtis-Hedlund-Lyndon Theorem \cite{H}:
any topological conjugacy  or semi-conjugacy $F:(X,\sigma) \to (Y,\sigma)$
between symbolic systems is given by a {\bf local rule} 
$f : \mathcal S_X^{m+1+a} \to \mathcal S_Y$,
where for every $x \in X$,
$(F(x))_i = f(x_{i-m},\dots,x_{i+a})$, 
$-\infty < i < \infty$. 
Here $m \ge 0$ is called the {\bf memory} 
and $a \ge 0$ the {\bf anticipation} of~$F$.
In this case $F$ is called an
$(m+1+a)$-block map.  
Note that by adding superfluous variables an $r$-block map is also an $s$-block map for all $s \ge r$.
The powers of $f$ are defined so that for every $n \ge 2$,  $f^n$ is a local rule for~$F^n$.  For example, if $f$ is a 3-block map, then $f^2$ is the 5-block map defined by
$$ f^2(x_1,\ldots,x_5) := f(f(x_1,x_2,x_3),f(x_2,x_3,x_4),f(x_3,x_4,x_5)).$$

A {\bf substitution of constant length} $L \ge 2$ is a 
mapping $\theta : \mathcal S \to \mathcal S^L$
from a finite set of symbols $\mathcal S$ to the set of 
$L$-blocks, i.e. blocks of length~$L$,
 with entries from $\mathcal S$.
The classic example is the {\bf Morse} substitution:
$0 \mapsto 01, 1 \mapsto 10$.
A substitution $\theta$ can be extended to a mapping of finite blocks
by concatenation.  
For example, $\theta(st) := \theta(s)\theta(t)$.
The powers of $\theta$ are defined in the obvious way.
For example, if $\theta(s) = tuv$,
then $\theta^2(s) := \theta(t)\theta(u)\theta(v)$.

A symbolic minimal system $(X,\sigma)$ is called a
{\bf substitution minimal system} if and only if
it is can be {\bf generated} by a
substitution~$\theta$, 
i.e. for every $s \in \mathcal S_X$,
the alphabet of~$X$, and for every  $n \ge 1$,
$\theta^n(s)$ appears in~$X$.
The substitution~$\theta$ defined on alphabet~$\mathcal S$
is called {\bf primitive}
if and only if there exists $n \ge 1$ such that
for every $s,t \in \mathcal S$,
$s$ appears in $\theta^n(t)$.

Lemmas 1-3 below hold for all substitutions, not just substitutions of constant  length.  Lemma~1 shows that when studying substitution minimal systems, there is no loss in generality in assuming that the generating substitution is primitive.

\begin{lemma1}\cite{DMK}\cite[Prop. 5.5]{Q}
Every primitive substitution generates a unique
substitution minimal system.
Conversely, every substitution minimal system can be generated by a primitive substitution.
\end{lemma1}

An important property of substitution minimal systems generated by primitive substitutions is
{\bf recognizability}, also called  
{\bf unique decipherability}, defined in the 
following lemma.

\begin{lemma2}\cite[Th\'{e}or\`{e}mes 1 et 2]{M}
Let $\theta$ be a primitive substitution %of constant length 
that generates an infinite symbolic minimal system $(X_\theta,\sigma)$.
Then
every point
in~$X_\theta$ can be written as a concatenation of
the blocks $\theta(s)$ in exactly one way.
\end{lemma2}

Lemma 2 is not true for  primitive substitutions that generate finite systems.  For example, $0 \mapsto 010$, $1 \mapsto 101.$

The unique substitution minimal system generated by
~$\theta$ is denoted ~$(X_\theta,\sigma)$.
A basic fact about substitution minimal systems is

\begin{lemma3}\cite[Prop. 5.4]{Q}
Let $(X_\theta,\sigma)$ be the substitution minimal system generated 
by~$\theta$.  Then for every $n \ge 1$,
$X_\theta = X_{\theta^n}$.
\end{lemma3}

Lemma 4 shows that we can assume, whenever it is useful, that a substitution~$\theta$ 
is {\bf one-to-one},
i.e. if $s \ne t$, then  $\theta(s) \ne \theta(t)$.

\begin{lemma4}\cite[Prop. 2.3]{BDM}
For any primitive, constant length substitution~$\theta$ such that 
$(X_\theta,\sigma)$ is infinite, there is a primitive, one-to-one substitution 
$\zeta$ of the same constant length such that
$(X_\theta,\sigma)$ and  $(X_\zeta,\sigma)$
are topologically conjugate.
\end{lemma4}

For  references on substitutions, see \cite{F} for arbitrary substitutions and \cite{G} for constant length substitutions.

\bigskip

\section*{3. Conjugacy to a given substitution minimal system}

In this section we find in Theorem 1 necessary and sufficient conditions for a symbolic minimal system to be topologically conjugate to a given constant length substitution minimal system.
For a subclass of substitutions, including the  
Morse substitution, we find  in Theorem~2 a result with a simpler statement and a simpler proof than for Theorem~1.

\begin{theorem1}  Let $\theta$ be a primitive, one-to-one substitution of constant length $L \ge 2$ such that $X_{\theta}$ is infinite, and let $(Y, \sigma)$ be a symbolic minimal system.  Then $(Y, \sigma)$ is topologically conjugate to $(X_{\theta}, \sigma)$ if and only if 

\begin{enumerate}
\item there exist $N \ge 1$ and 
a collection $\mathcal B$ of $L^N$-blocks such that
every point in $Y$ can be written  as a concatenation of
blocks in $\mathcal B$ 
in exactly one way, 

\item there exists a 2-block semiconjugacy  
$G: (X_\theta,\sigma) \to (Y_0,\sigma^{L^N})$,
where $Y_0$ is the set of points in~$Y$ such that the 
blocks in $\mathcal B$
start at multiples of~$L^N$,
and

\item if $stu$ and $s't'u'$ are 3-blocks appearing in 
~$X_\theta$ and $t \ne t'$, then $g(stu) \ne g(s't'u')$,
where $g$ is a local rule of~$G$.
\end{enumerate}

\end{theorem1}

\begin{proof} (\emph{only if}) Let $F: (X_\theta,\sigma) \to (Y,\sigma)$ be a topological conjugacy.  By composing  with a power of the shift, we may assume that $F$ has no memory.  Let $f$ be a local rule of~$F$.
Choose  $N$ so  that 
 $L^N$ is greater than the anticipation of~$F$ and
\medskip
\newline 
($*$) for all 3-blocks
$stu$ and $s't'u'$  appearing in $X_\theta$ with $t \ne t'$,
$f(\theta^N(s) \theta^N(t) \theta^N(u)) \ne
f(\theta^N(s') \theta^N(t') \theta^N(u'))$.
\medskip
\newline
To see that  $N$ can be chosen so that ($*$) holds,
suppose not.
Then equality holds for some $stu$ and $s't'u'$
with $t \ne t'$
and infinitely many $n$.
For every ~$n$ and every $s,t,u$,
$\theta^n(s) \theta^n(t) \theta^n(u)$
can be extended to a doubly infinite sequence,
that  we also call
$\theta^n(s) \theta^n(t) \theta^n(u)$,
with the zeroth coordinate coming in a place of
disagreement between $\theta^n(t)$ and ~$\theta^n(t').$
Using the fact that
$\theta$ and hence every $\theta^n$ is one-to-one, 
there is a subsequence $(\tilde{n})$ of $(n)$
such that both
$(\theta^{\tilde{n}}(s) \theta^{\tilde{n}}(t) \theta^{\tilde{n}}(u))$ and
$(\theta^{\tilde{n}}(s') \theta^{\tilde{n}}(t') \theta^{\tilde{n}}(u'))$
converge, say to $x \ne x'$.
Then $F(x) = F(x')$,
contradicting $F$ being one-to-one.

By adding superfluous variables if necessary, we may assume that the anticipation of ~$F$  is exactly ~$L^N$.  
Define $$\mathcal B := \{f(\theta^N(st)): st \text{ is a 2-block appearing in } X_\theta\}.$$

We show that condition ~(1) holds.  Since $\theta$
and hence $\theta^N$ are primitive, every point in $X_{\theta^N} = X_\theta$ can be written as a concatenation of blocks  $\theta^N(s)$ in exactly one way \cite{M}.
Then, since $F$ is a topological conjugacy,
every point in $Y$ can be written as a concatenation of
blocks in $\mathcal B$
in exactly one way.

To show that condition~(2) holds, let $G : (X_\theta,\sigma) \to (Y_0,\sigma^{L^N})$, where $Y_0$ 
is the set of points in~$Y$ such that the blocks in~$\mathcal B$ start at multiples of~$L^N$, be the 2-block semiconjugacy with local rule
$g := f \circ \theta^N$. 

Condition~(3) follows from ($*$).

\medskip

(\emph{if}) It follows from~(3) that the topological semiconjugacy $G$ is one-to-one and hence a topological conjugacy.  

Let $X_{\theta,0}$ be the set of points in $X_\theta$
such that the blocks $\theta^N(s)$ start at multiples
of~$L^N$.  Then $(X_\theta,\sigma)$ is topologically conjugate to
$(X_{\theta,0},\sigma^{L^N})$ via the map
$s \mapsto \theta^N(s)$.

Standard arguments (see, e.g.,  the proof of the Morse Dynamical Characterization Theorem in \cite{CKL}) show that
$(X_\theta,\sigma)$ is topologically conjugate to
$(Y,\sigma)$.
\end{proof}

Now we consider a special case, substitutions  $\theta$ of constant length at least three (see Remark below) for which for all $s \ne t$,
$\theta(s)$ and $\theta(t)$ disagree in some entry other than the first or last.  This class contains the square of the Morse substitution.

For such substitutions, condition $(*)$ in the proof of Theorem~1 can be replaced by the stronger condition $(**)$ below 
and  we have the following improvement of Theorem~1 and generalization of the Morse Dynamical Characterization Theorem~\cite{CKL}.

\begin{theorem2}
Let $\theta$ be a primitive, one-to-one substitution of constant length $L \ge 3$ such that $X_{\theta}$ is infinite, and let $(Y, \sigma)$ be a symbolic minimal system.
Suppose also that
any two substitution blocks  $\theta(s)$ and  $\theta(t)$ with $s \ne t$ disagree somewhere other than in the first or last entry.

Then $(Y, \sigma)$ is topologically conjugate to $(X_{\theta}, \sigma)$ if and only if
there exist  $N \ge 1$, $a \ge 0$, a collection $\mathcal B$ of 
$(L^N - a)$-blocks
that are in one-to-one correspondence with the symbols in~$X_\theta$, 
and a collection $\mathcal B'$ of $a$-blocks  
such that
\begin{enumerate}
\item
every point in $Y$ can be written as a concatenation of alternating  
blocks in~$\mathcal B$ and blocks in~$\mathcal B'$
in exactly one way,
\item the blocks in ~$\mathcal B'$ are determined by their nearest
neighbors in~$\mathcal B$,
\item every second block in %the 
a concatenation~(1),
thought of as an infinite bilateral sequence 
with letters from ~$\mathcal B$,
``mirrors'' a point in~$X_\theta$
via the  one-to-one correspondence above.
\end{enumerate}
\end{theorem2}

\begin{remark}
We require $L \ge 3$ because condition~(3) is vaccuous for $L=2$.  However, this requirement is harmless, for in this case look at $\theta^2$ rather than~$\theta$.
\end{remark}

The proof of the \emph{only if} direction proceeds as in the proof of Theorem~1, with
$F:(X,\sigma) \to (Y,\sigma)$ being a topological conjugacy with no memory and anticipation~$a$ and 
condition~$(**)$ below taking the place of condition~$(*)$. 
\medskip
\newline 
$(**)$ there exists $N$ such that $L^N$ is greater than
the anticipation of ~$F$ and such that for all symbols
$t \ne t'$ appearing in~$X_\theta$,
$f(\theta^N(t)) \ne f(\theta^N(t'))$.
\medskip
\newline 
Let $\mathcal B := \{f(\theta^N(s))\}$ and let
$\mathcal B'$ be the $a$-blocks appearing in points of~$Y$
between consecutive blocks in~$\mathcal B$ .

To see that condition (2) holds, note that every
block in~$\mathcal B'$ appears in some
$f(\theta^N(st))$.

  The proof of the  \emph{if} direction 
is much the same as it is in the proof of Theorem~1.

\bigskip

\section*{4. Number of symbols}

In this section we find a relation between the number of symbols in constant length substitutions that generate topologically conjugate systems.

\begin{theorem3}
Let $\theta$ be a primitive, constant length substitution such that $X_\theta$ is infinite.  Then the number of symbols in a primitive, one-to-one constant length substitution that generates a substitution minimal system topologically  conjugate to $(X_\theta,\sigma)$ is at most the number of 3-blocks appearing in~$X_\theta$.
\end{theorem3}

\begin{proof}
By \cite{CDL} we may assume that the lengths of $\theta$ and~$\zeta$ are the same.  
Let $F : (X_\theta,\sigma) \to (X_\zeta,\sigma)$
be a topological conjugacy with no memory and local rule~$f$.
Then, with notation as in the proof of Theorem~1,
for every symbol $s'$ appearing in~$X_\zeta$,
$\zeta^N(s')$ is a subblock of some $B_1B_2$,
where $B_1,B_2 \in \mathcal B$.
But $B_1B_2 = f(\theta^N(stu))$ for some 3-block $stu$
appearing in~$X_\theta$.
Then we have
$$ \#s' = \#\zeta^N(s')
		      		       \le \# f(\theta^N(stu))
		       \le \# stu,
$$
the equality because $\zeta$ is one-to-one, and the first inequality because $\theta$ and $\zeta$ are uniquely decipherable.
\end{proof}

The following six-symbol example shows that the bound is %sharp 
attained
for the Morse substitution.  It is one of the ``3-block presentations"  of the Morse substitution and so generates a system topologically conjugate to the Morse system (see \cite{CDK}).

\medskip

\begin{align}
\zeta(001) = (101)(011)  \notag \\
\zeta(010) = (110)(100)  \notag \\
\zeta(011) = (110)(101)  \notag \\
\zeta(100) = (001)(010)  \notag \\
\zeta(101) = (001)(011)  \notag \\
\zeta(110) = (010)(100)  \notag
\end{align}

With $\theta$ the Morse substitution
$0 \mapsto 01$, $1 \mapsto 10$,
$\zeta(stu) := (s_2 t_1 t_2)(t_1 t_2 u_1)$,
where $\theta(s) = s_1 s_2$, etc.

It follows from \cite[Toeplitz Corollary]{CKL} that the bound (five) given by Theorem~3 for the Toeplitz substitution 
$0 \mapsto 01, 1 \mapsto 00$ cannot be attained.

\end{document}